\documentclass[12pt,twoside]{article}
\usepackage{amssymb,mathrsfs,amsmath}

\long\def\@makefntext#1{\noindent #1}
\newskip\tabcentering \tabcentering=1000pt plus 1000pt minus 1000pt
\def\REF#1{\par\hangindent\parindent\indent\llap{#1\enspace}\ignorespaces} %reference format
\def\MCH#1#2{\setbox0=\hbox{\raise#1\hbox{#2}}\smash{\box0}}% move char

\def\@evenfoot{}\def\@oddfoot{}

\def\@evenhead{\hbox to\textwidth{\footnotesize\hfill\rm\thepage
}} % authors name

\def\@oddhead{\hbox to \textwidth{\footnotesize \hfill\rm\thepage}}% abbreviate title

\def\sec#1{\vskip 3mm\leftline{\bf #1}\vskip 1mm}

\def\bc{\begin{center}}
\def\ec{\end{center}}
\def\no{\noindent}
\def\hang{\hangindent\parindent}
\def\textindent#1{\indent\llap{\qquad #1\ \ \enspace}\ignorespaces}
\def\ref{\par\hang\textindent}

\begin{document}
\abovedisplayskip=6pt plus 1pt minus 1pt \belowdisplayskip=6pt
plus 1pt minus 1pt
%-------------------  First Head  -----------------------------------------
\thispagestyle{empty} \vspace*{-1.0truecm} \noindent
%\parbox[b]{6truecm}{\footnotesize\baselineskip=11pt\noindent Acta Mathematica
%Sinica, English Series\\

%Http://www.ActaMath.com} \hfill
%\parbox[t]{6truecm}{\vskip -1.7cm \hfill\includegraphics{actmark.eps}}
%===================Text=============================================
\vskip 8mm

\bc{\large\bf The thickness of amalgamations of graphs$^\ast$

\footnotetext{\footnotesize $^\ast$Supported by NNSF of China
under Grant
No.11001196 and No.11126167\\
$^\dag$Corresponding author.\
 {\it E-mail address:}
yanyang@tju.edu.cn (Y.Yang)}} \ec

\vskip 3mm \bc{\bf Yan Yang$^\dag$ \ and \ Xiangheng Kong\vskip 1mm {\small\it
Department of Mathematics, Tianjin University, Tianjin {\rm
300072,} P.R.China }}\ec \vskip 1 mm

\noindent{\small {\small\bf Abstract} \ \ The thickness $\theta(G)$ of a graph $G$ is the minimum number of planar spanning subgraphs  into
which the graph $G$ can be decomposed. As a topological invariant of a graph, it is a measurement of the closeness to planarity of a graph, and it also has important applications to VLSI design. In this paper, the thickness of graphs that are obtained by vertex-amalgamation and
bar-amalgamation of any two graphs whose thicknesses are known are obtained, respectively. And the lower and upper bounds for the thickness of graphs that are obtained by edge-amalgamation and 2-vertex-amalgamation of any two graphs whose thicknesses are known are also derived, respectively.

\vspace{2mm}

\no{\small\bf Keywords} \ \ thickness; amalgamation; genus. \vspace{2mm}

 \vspace{5mm}

 \sec{\large 1\quad Introduction}\vspace{1mm}

\no Let $G$ be a graph with  vertex set $V(G)$ and edge set $E(G)$. A graph is said to be {\it planar} if it can be  drawn on the plane so that no two edges cross (i.e., its edges meet only at their common ends); otherwise, non-planar. As for a non-planar graph, there are some measurements of the closeness to planarity of a graph, such as thickness, genus, crossing number etc.

Suppose $G_1, G_2, \ldots, G_k$ are spanning subgraphs of $G$, if $$E(G_1)\cup E(G_2)\cup\cdots \cup E(G_k)=E(G) ~\mbox{and} ~ E(G_i)\cap E(G_j)=\emptyset,  (i\neq j,~i,j=1,2,\ldots,k),$$ then $\{G_1, G_2, \ldots, G_k\}$ is a decomposition of $G$. Furthermore, if $G_1, G_2, \ldots, G_k$ are all planar graph, then $\{G_1, G_2, \ldots, G_k\}$ is a planar decomposition of $G$. The minimum number of planar spanning subgraphs into
which a graph $G$ can be decomposed is called the {\it thickness} of $G$, denoted by $\theta(G)$.

As a  topological invariant of a graph, thickness is an important research object in topological graph theory. And
it also has important applications to VLSI design[1]. But the results about thickness are few,  compared with other topological invariants, e.g., genus, crossing number. The only types of graphs whose thicknesses have been obtained are complete graphs[6], complete bipartite graphs[7] and hypercubes[13]. Since determining the thickness of graphs is NP-hard[14], it is very difficult to get the exact number of thickness for arbitrary graphs, people study lower and upper bounds for the thickness of a graph[10,12] and introduce  heuristic algorithms to approximate it[9,15]. And some relations between thickness and other topological invariants, such as genus, are also established[2].

In this paper,  the thickness of graphs that are obtained by vertex-amalgamation and bar-amalgamation of any two graphs whose thicknesses are known are obtained, respectively. And the lower and upper bounds for the thickness of graphs that are obtained by edge-amalgamation
and 2-vertex-amalgamation of any two graphs whose thicknesses are known are also derived, respectively.

Graphs in this paper are simple graphs. For the undefined terminologies see [5].

\vspace{5mm}

 \sec{\large 2\quad Thickness of graph amalgamations}\vspace{1mm}

\no The union of graphs $G_1$ and $G_2$ is the graph $G_1\cup G_2$
with vertex set $V (G_1) \cup V (G_2)$ and edge set $E(G_1) \cup E(G_2)$.
The intersection $G_1 \cap G_2$ of $G_1$ and $G_2$ is defined analogously.

Let $G_1$ and $G_2$ be subgraphs of a graph $G$. If $G=G_1\cup G_2$ and $G_1 \cap G_2=\{v\}$ (a vertex of $G$),
then we say that $G$ is the {\it vertex-amalgamation} of $G_1$ and $G_2$ at vertex $v$, denoted  $G=G_1\vee^{1}_{\{v\}}G_2$.
 If $G=G_1\cup G_2$ and $G_1\cap G_2=\{v,u\}$ (two distinct vertices of $G$), then we say that $G$ is the {\it 2-vertex-amalgamation} of $G_1$ and $G_2$ at vertices $v$ and $u$, denoted  $G=G_1\vee^{1}_{\{v,u\}}G_2$. If $G=G_1\cup G_2$ and $G_1 \cap G_2=\{e\}$ (an edge of $G$),
then we say that $G$ is the {\it edge-amalgamation} of $G_1$ and $G_2$ on edge $e$, denoted  $G=G_1\vee^{2}_{\{e\}}G_2$.

Let $G$ and $H$ be two disjoint graphs, the {\it bar-amalgamation} of $G$ and $H$ is obtained by running an new edge between a vertex of $G$ and a vertex of $H$.

The four kinds of amalgamations defined above are important operations on graphs, by these amalgamations, one can get larger graphs
(i.e., the graph with larger order) synthesized from small graphs. It is a general method to study problems in graph theory  by operations on graphs.
For example, we will list some results about genus of graph amalgamations in the following.

The {\it genus} of a graph $G$ is the minimum integer $k$ such that $G$ can be embedded on the orientable surface of genus $k$, denoted by $\gamma(G)$. A graph $G$ is planar if and only if $\gamma(G)=0.$

\vskip1mm

\no{\bf Lemma 2.1}[4]\quad{\it  If $G$ is the vertex-amalgamation of $G_1$ and $G_2$, then $$\gamma(G)=\gamma(G_1)+\gamma(G_2).$$}
\no{\bf Lemma 2.2}[8]\quad{\it If $G$ is the bar-amalgamation of $G_1$ and $G_2$, then $$\gamma(G)=\gamma(G_1)+\gamma(G_2).$$}
\no{\bf Lemma 2.3}[3]\quad{\it If $G$ is the edge-amalgamation of $G_1$ and $G_2$, then $$\gamma(G)\leq\gamma(G_1)+\gamma(G_2).$$}
\no{\bf Lemma 2.4}[11]\quad{\it  If $G$ is the 2-vertex-amalgamation of $G_1$ and $G_2$, then
 $$\gamma(G_1)+\gamma(G_2)-1\leq\gamma(G)\leq\gamma(G_1)+\gamma(G_2)+1.$$}

In [2], a relation between genus and thickness of a graph is given.
\vskip1mm
\no{\bf Lemma 2.5}[2]\quad{\it Let $G$ be a simple graph, if $\gamma(G)=1$, then $\theta(G)=2$.}\vskip2.5mm

In the following, some theorems about the thickness of vertex-amalgamation, bar-amalgamation, edge-amalgamation
 and  2-vertex-amalgamation of graphs are obtained.
\vskip1.5mm
 \vskip1mm
\vskip1mm
 \no{\bf Theorem 2.1}\quad{\it If $G$ is the vertex-amalgamation of $G_1$ and $G_2$, $\theta(G_1)=n_1$ and $\theta(G_2)=n_2$, then $$\theta(G)= \max\{n_1, n_2\}.$$}
{\it Proof}\quad Without loss of generality, one can assume that $n_1\geq n_2$ and $G_1 \cap G_2=\{v\}$  (a vertex of $G$).
Suppose that
$\{G_{11}, G_{12}, \ldots, G_{1n_1}\}$ is a planar decomposition of $G_1$ and
$\{G_{21}, G_{22}, \ldots, G_{2n_1}\}$ is a planar decomposition of $G_2$. From Lemma 2.1,
$$\gamma(G_{1i}\vee^{1}_{\{v\}} G_{2i})=\gamma(G_{1i})+\gamma(G_{2i})=0+0=0,~~~~~ 1\leq i\leq n_1.$$
Hence $\{~(G_{11}\vee^{1}_{\{v\}} G_{21}), ~ (G_{12}\vee^{1}_{\{v\}} G_{22}), \ldots, (G_{1n_1}\vee^{1}_{\{v\}} G_{2n_1})~\}$
is a planar decomposition of $G$, ~$\theta(G)\leq n_1$. On the other hand, $G=G_{1}\vee^{1}_{\{v\}} G_{2}$, ~$G_1$ is a subgraph of $G$ and $\theta(G_1)=n_1$, so $\theta(G)\geq n_1$. Summarizing the above, $\theta(G)=n_1$, the theorem follows.\hspace{\fill}$\Box$

\vskip1.5mm
 \vskip1mm
\vskip1mm
 \no{\bf Theorem 2.2}\quad{\it  If $G$ is the bar-amalgamation of $G_1$ and $G_2$, $\theta(G_1)=n_1$ and  $\theta(G_2)=n_2$, then
 $$\theta(G)= \max\{n_1, n_2\}.$$}
{\it Proof}\quad Suppose that $n_1\geq n_2$ and edge $e$ is the new edge between  $G_1$ and  $G_2$.
Let  $\{G_{11}, G_{12}, \ldots, G_{1n_1}\}$ be a planar decomposition of $G_1$ and
$\{G_{21}, G_{22}, \ldots, G_{2n_1}\}$ be a planar decomposition of $G_2$. $G_{11}\cup G_{21}\cup {e}$ is
the bar-amalgamation of $G_{11}$ and $G_{21}$, from Lemma 2.2, the genus of $G_{11}\cup G_{21}\cup {e}$ is zero, that is to say,
 $G_{11}\cup G_{21}\cup {e}$ is a planar graph. Hence $\{G_{11}\cup G_{21}\cup {e}, ~G_{12}\cup G_{22} , \ldots, G_{1n_1}\cup G_{2n_1}\}$
 is a planar decomposition of $G$, ~$\theta(G)\leq n_1$. And $G=G_{1}\cup G_{2}\cup {e}$, $\theta(G_1)=n_1$, so $\theta(G)\geq n_1$.
Summarizing the above, $\theta(G)=n_1$, the theorem is obtained.\hspace{\fill}$\Box$

 \vskip1.5mm
 \vskip1mm

\vskip1mm
 \no{\bf Theorem 2.3}\quad{\it  If $G$ is the edge-amalgamation of $G_1$ and $G_2$, $\theta(G_1)=n_1$ and  $\theta(G_2)=n_2$, then}
 $$ \max\{n_1, n_2\}\leq\theta(G)\leq n_1+n_2-1.$$
 \vskip1mm
{\it Proof}\quad Suppose that $n_1\geq n_2$ and $G_1 \cap G_2=\{e\}$ (an edge of $G$), the two end vertices of $e$ are $v$ and $u$.
Let  $\{G_{11}, G_{12}, \ldots, G_{1n_1}\}$ be a planar decomposition of $G_1$ and
$\{G_{21}, G_{22}, \ldots, G_{2n_2}\}$ be a planar decomposition of $G_2$. Without loss of generality, one can assume that
the edge $e$ contains in $G_{11}$ and $G_{21}$.  From Lemma 2.3, $G_{11}\vee^{2}_{\{e\}}G_{21}$ is a planar graph, and from Lemma
2.4 and 2.5,
$$\gamma(G_{1i}\vee^{1}_{\{v,u\}}G_{2i})\leq 1 ~~\mbox{and}~~\theta(G_{1i}\vee^{1}_{\{v,u\}}G_{2i})\leq 2, ~~~2\leq i\leq n_2.$$ Let $\{G^i,\widetilde{G^i}\} $ be a  planar decomposition of $G_{1i}\vee^{1}_{\{v,u\}}G_{2i}$, for $2\leq i\leq n_2$, then
$$\{~G_{11}\vee^{2}_{\{e\}}G_{21},~ G^2,~\widetilde{G^2}, \ldots, G^{n_2},~\widetilde{G^{n_2}},~ G_{1n_2+1},\ldots, G_{1n_1}~\} $$
 is  a planar decomposition of $G$, therefor $\theta(G)\leq 1+2(n_2-1)+n_1-n_2=n_1+n_2-1.$ And $G=G_{1}\vee^{2}_{\{e\}}G_{2}$, $\theta(G_1)=n_1$, so $\theta(G)\geq n_1$. Summarizing the above, the theorem follows.\hspace{\fill}$\Box$

\vskip1.5mm
 \vskip1mm

\vskip1mm
 \no{\bf Theorem 2.4}\quad{\it  If $G$ is the 2-vertex-amalgamation of $G_1$ and $G_2$, $\theta(G_1)=n_1$ and  $\theta(G_2)=n_2$, then}
 $$ \max\{n_1, n_2\}\leq\theta(G)\leq n_1+n_2.$$
 \vskip1mm
 \vskip1mm
{\it Proof}\quad With a similar argument to the proof of Theorem 2.3, one can obtain the theorem easily.
\hspace{\fill}$\Box$

 \vskip1.5mm
 \vskip1mm

\no \vskip0.3in \no {\bf References} \vskip0.1in

\footnotesize \

\REF{[1]} A.Aggarwal, M.Klawe, and P.Shor, Multilayer grid embeddings for VLSI, Algorithmica, 6(1991) 129-151го

\REF{[2]} K.Asano, On the genus and thickness of graphs, J. Comb. Theory (B), 43 (1987) 187-192го

\REF{[3]} S.Alpert, The genera of amalgamations of graph, Trans. Amer. Math. Soc., 178(1973) 1-39.

\REF{[4]} J.Battle, F.Harary, Y.Kodama and J.W.T.Youngs, Additivity of the genus of a graph, Bull. Amer. Math. Soc., 68 (1962) 565-568.

\REF{[5]} J.A.Bondy and U.S.R.Murty, Graph Theory, GTM 244, Springer, 2008го

\REF{[6]} L.W.Beineke and F.Harary, The thickness of the complete graph, Canad. J. Math., 17 (1965) 850-859го

\REF{[7]} L.W.Beineke, F.Harary, and J. W.Moon, On the thickness of the complete bipartite graph, Proc. Cambridge. Philo. Soc., 60 (1964) 1-5го

\REF{[8]} J.E.Chen, S.P.Kanchi and A.Kanevsky, A note on approximating graph genus,  Inform. Proc. Letters, 61(1997) 317-322.

\REF{[9]} R.J.Cimikowski, On heuristics for determining the thickness of a graph, Inform. Sci., 85 (1995) 87-98го

\REF{[10]} A.M.Dean, J.P.Hutchinson, and E.R.Scheinerman, On the thickness and arboricity of a graph, J. Comb. Theory (B), 52(1991) 147-151го

\REF{[11]} R.Decker, H.Glover and J.P.Huneke,  The genus of the 2-amalgamations of graphs, J. Graph Theory, 5(1981) 95-102.

\REF{[12]} J.H.Halton, On the thickness of graphs of given degree, Inform. Sci., 54 (1991) 219-238го

\REF{[13]} M. Kleinert,  Die Dicke des n-dimensionale W$\ddot{\mbox{u}}$rfel-Graphen, J. Comb. Theory, 3 (1967) 10-15.

\REF{[14]} A.Mansfield, Determining the thickness of graphs is NP-hard, Math. Proc. Cambridge Philo. Soc., 93: 9 (1983) 9-23го

\REF{[15]} T.Poranen, A simulated annealing algorithm for determining the thickness of a graph, Inform. Sci., 172 (2005) 155-172го

\end{document}